\documentclass[12pt]{article}
\usepackage[a4paper, total={7in, 10in}]{geometry}
\usepackage[parfill]{parskip}
\usepackage{physics, tensor, float, subcaption}
\usepackage{graphicx}
\graphicspath{ {Plots/} }
\usepackage{jhep-mod}
\usepackage{bm}
\usepackage{soul}
\usepackage{amssymb,amsmath,amsthm}
\usepackage{mathrsfs}
\usepackage[utf8]{inputenc}
\usepackage{enumerate}
\usepackage{bigints}
\usepackage{xcolor}
\usepackage{appendix}
\usepackage{graphicx}
\usepackage{float}
\usepackage{tikz}
\usepackage{setspace}
\usepackage{cancel}
\usepackage{doi}
\definecolor{purple}{rgb}{1,0,1}
\definecolor{lime}{HTML}{A6CE39} 

\newcommand{\blue}[1]{{\color{blue} #1}}

\definecolor{lime}{HTML}{A6CE39}
\newcommand{\orcidicon}{%
	\begin{tikzpicture}
	\draw[lime, fill=lime] (0,0) 
		circle [radius=0.16] 
		node[white] {{\fontfamily{qag}\selectfont \tiny ID}};
	\draw[white, fill=white] (-0.0625,0.095) 
		circle [radius=0.007];
	\end{tikzpicture}
	\hspace{-5mm}
}
\newcommand\orcidMatt{{\href{https://orcid.org/0000-0003-1088-6485}{\orcidicon}}}

\renewcommand{\O}{\mathcal{O}}
\begin{document}

\title{\null\vspace{-50pt}
\leftline{The $n^{th}$ prime exponentially}
}

\author{
\Large
Matt Visser\!\orcidMatt\!
}
\affiliation{School of Mathematics and Statistics, Victoria University of Wellington, \\
\null\qquad PO Box 600, Wellington 6140, New Zealand.}
\emailAdd{matt.visser@sms.vuw.ac.nz}
\def\theta{\vartheta}
\def\O{{\mathcal{O}}}
\def\Li{{\mathrm{Li}}}

\abstract{
\vspace{1em}

From  known effective bounds on the prime counting function of the form
\[
|\pi(x)-\Li(x)| <  a  \;x \;(\ln x)^{b} \; \exp\left(-{c}\; \sqrt{\ln x}\right); 
\qquad (x \geq x_0);
\]
it is possible to establish exponentially tight effective upper and lower bounds on the prime number theorem:
For $x \geq x_*$ where  $x_*\leq \max\{x_0,17\}$ we have:
\[
{\Li(x) \over 1+a\; (\ln x)^{b+1} \; \exp\left(-c\; \sqrt{\ln x}\right)} < \pi(x) < {\Li(x) \over 1-a \;(\ln x)^{b+1} \; \exp\left(-c\; \sqrt{\ln x}\right)}.
\]
Furthermore, it is possible to establish exponentially tight effective upper and lower bounds on 
the location of the $n^{th}$ prime. Specifically:
\[
p_n < \Li^{-1} \left( n \left[1+  a   \;(\ln[n\ln n])^{b+1} \; \exp\left(-{c}\; \sqrt{\ln[n\ln n]}\right)\right] \right);
\qquad (n\geq n_*).
\]
\[
p_n > \Li^{-1} \left( n \left[1-  a   \;(\ln[n\ln n])^{b+1} \; \exp\left(-{c}\; \sqrt{\ln[n\ln n]}\right)\right] \right);
\qquad (n\geq n_*).
\]
Here the range of validity is explicitly bounded by some $n_*$ satisfying
\[
n_* \leq \max\left\{\pi(x_0),\pi(17), \pi\left( (1+e^{-1}) \exp\left( \left[2(b+1)\over c\right]^2\right)\right) \right\}.
\]
Many other fully explicit bounds along these lines can easily be developed. 

\bigskip

\bigskip
\noindent
{\sc Date:} 20 April 2025; \LaTeX-ed \today


}

\maketitle
\def\tr{{\mathrm{tr}}}
\def\diag{{\mathrm{diag}}}
\def\cof{{\mathrm{cof}}}
\def\pdet{{\mathrm{pdet}}}
\def\d{{\mathrm{d}}}
\parindent0pt
\parskip7pt
\def\Kerr{{\scriptscriptstyle{\mathrm{Kerr}}}}
\def\eos{{\scriptscriptstyle{\mathrm{eos}}}}
\section{Introduction: Known bounds on $|\pi(x)-\Li(x)|$}

Work over the last decade or so has developed a number of fully effective bounds on the prime counting function $\pi(x)$ of the general form~\cite{Trudgian,Johnston-Yang,Fiori-et-al}:
\begin{equation}
|\pi(x)-\Li(x)| <  a  \;x \;(\ln x)^{b} \; \exp\left(-{c}\; \sqrt{\ln x}\right); 
\qquad (x \geq x_0).
\end{equation}

\smallskip
For some widely applicable effective bounds of this type see Table I below.\\
(A straightforward elementary numerical computation is required to determine the numerical coefficients in the Trudgian~\cite{Trudgian} bounds.)

\begin{table}[!h]
\caption{Some widely applicable effective bounds on $\pi(x)$.}\smallskip
\begin{center}
\begin{tabular}{||c|c|c|c||c||c||}
\hline
\hline
$a$ & $b$ & $c$ & $x_0$ & Source & Notes\\
\hline
\hline
0.4394& --3/4 &0.32115 & 59 & Trudgian~\cite{Trudgian} & Eq (7) \\ 
\hline
0.2795& --3/4 &0.3936 & 229 & Trudgian~\cite{Trudgian} & Th 2 \\ 
\hline
\hline
9.2211 & 1/2 & 0.8476 & 2 & Fiori--Kadiri--Swidinsky~\cite{Fiori-et-al}& Eq (3), Eq (43)\\ %
\hline
\hline
9.59 & 0.515 & 0.8274 & 2 & Johnston--Yang~\cite{Johnston-Yang}& Eq (1.6)\\
\hline
\hline 
\end{tabular}
\end{center}
\end{table}
 
\clearpage
\section{Some elementary variants on these bounds}

\enlargethispage{40pt}
There are various strategies for relaxing these bounds to make them more tractable. The most simple technique is to simply increase either of the coefficients $\{a,b\}$, or decrease the coefficient $c$, in the Trudgian bounds so as to  increase the domain of validity (decrease $x_0$). 
(The FKS and JY bounds are already maximal in this regard.)

A more subtle technique (based on a variant of the discussion in reference~\cite{dlVP}) is as follows: Let $\tilde b\leq b$ and $\tilde c < c$.
Write
\begin{eqnarray}
&& a \;x \;(\ln x)^b \; \exp\left(-c\; \sqrt{\ln x}\right) \nonumber\\
&& =
a \;x \left\{ (\ln x)^{b-\tilde b} \; \exp\left(-[c-\tilde c]\; \sqrt{\ln x}\right) \right\} \;(\ln x)^{\tilde b} \; \exp\left(-\tilde c\; \sqrt{\ln x}\right).
\end{eqnarray}
Then the quantity in braces is bounded, reaching a global maximum at 
\begin{equation}
x_{max} = \exp\left( 4[b-\tilde b]^2\over[c-\tilde c]^2\right),
\end{equation}
where it takes on the value
\begin{eqnarray}
\left\{ (\ln x)^{b-\tilde b} \; \exp\left(-[c-\tilde c]\; \sqrt{\ln x}\right) \right\}_{max} 
&=&
\left( 4[b-\tilde b]^2\over[c-\tilde c]^2\right)^{b-\tilde b} \exp(-2[b-\tilde b])
\end{eqnarray}
Then defining
\begin{equation}
\tilde a = a\; \left( 4[b-\tilde b]^2\over[c-\tilde c]^2\right)^{b-\tilde b} \exp(-2[b-\tilde b])
= a\left( 2[b-\tilde b]\over e [c-\tilde c]\right)^{2(b-\tilde b)}.
\end{equation}
we see that
\begin{equation}
a \;x \;(\ln x)^b \; \exp\left(-c\; \sqrt{\ln x}\right) \leq 
\tilde a \;x \;(\ln x)^{\tilde b} \; \exp\left(-\tilde c\; \sqrt{\ln x}\right);
\qquad (x>1);
\end{equation}
with equality only at $x=x_{max}$.
Consequently any bound of the form
\begin{equation}
|\pi(x)-\Li(x)| <  a  \;x \;(\ln x)^{b} \; \exp\left(-{c}\; \sqrt{\ln x}\right); 
\qquad (x \geq x_0);
\end{equation}
implies bounds of the form
\begin{equation}
|\pi(x)-\Li(x)| <  \tilde a  \;x \;(\ln x)^{\tilde b} \; \exp\left(-{\tilde c}\; \sqrt{\ln x}\right); 
\qquad (x \geq x_0).
\end{equation}

Here the coefficients $\{\tilde a, \tilde b, \tilde c\}$ are as defined above. (Both $\tilde b< b$
and $\tilde c < c$ are free to be chosen, while $\tilde a > a$ is to be derived.) As long as $x_0$ is not too large, one can use explicit computations to find a larger domain of validity (some smaller $\tilde x_0$). 

Specifically there will be some $\tilde x_0 \leq x_0$ such that
\begin{equation}
|\pi(x)-\Li(x)| <  \tilde a  \;x \;(\ln x)^{\tilde b} \; \exp\left(-{\tilde c}\; \sqrt{\ln x}\right); 
\qquad (x \geq \tilde x_0; \; \tilde x_0 \leq x_0).
\end{equation}

In previous work, see reference~\cite{dlVP}, I had focussed on the case $\tilde b=0$ which is most appropriate to finding explicit bounds (for $\theta(x)$) of the de la Vall\'e Poussin form~\cite{Poussin}. Herein we shall soon see that for technical reasons when considering $\pi(x)$ it is more useful to consider the case $\tilde b = -1$. 
Table II below presents a number of derived widely applicable effective bounds of this general type. 

\begin{table}[!h]
\caption{Several derived widely applicable effective bounds on $\pi(x)$.\\
\centerline{(Input parameters on the left; output parameters on the right.)}\\
\centerline{(Note the bounds are closest to each other in the vicinity of $x_{max}$.)}\\
\centerline{$|\pi(x)-\Li(x)| <  \tilde a  \;x \;(\ln x)^{\tilde b} \; \exp\left(-{\tilde c}\; \sqrt{\ln x}\right).$}
}\smallskip
\begin{center}
\begin{tabular}{|||c|c|c|c|||c|c|c|c|c|||}
\hline
\hline
\vphantom{\Big |}
$a$ & $b$ & $c$ & $x_0$ & $\tilde a$ & $\tilde b$ & $\tilde c$ & $x_{max} $ & $\tilde x_0 $ \\
\hline
\hline
0.4394& --3/4 &0.32115 & 59 & 0.4680 & --7/8 & 1/4 & 203 931 & 41 \\ 
\hline
 0.4394& --3/4 &0.32115 & 59 & 0.4795 & --1 & 1/6 & 35 439 & 41\\
\hline
\hline
0.2795& --3/4 &0.3936 & 229 & 0.2804 &--5/6 & 1/3 & 2 097 & 227 \\ 
\hline
0.2795& --3/4 &0.3936 & 229 & 0.3164  & --1 & 1/4 & 184 165 & 223\\
\hline
\hline
9.2211 & 1/2 & 0.8476 & 2 & 9.7590 & 0 & 1/2 &3930 & 2\\ %
\hline
9.2211 & 1/2 & 0.8476 & 2 & 11.9026 & --1/2 & 1/5 & 13 874 & 2 \\ %
\hline 
9.2211 & 1/2 & 0.8476 & 2 & 29.6698 & --1 & 1/10 & 9 849 130 & 2 \\ %
\hline
\hline
9.59 & 0.515 & 0.8274 & 2 & 11.148 & 0 & 1/2 & 19 877 & 2\\
\hline
9.59 & 0.515 & 0.8274 & 2 & 13.659 & --1/2 & 1/5 & 35 206 & 2\\
\hline
9.59 & 0.515 & 0.8274 & 2 & 34.955 & --1 & 1/10 & 34 331 213 & 2\\
\hline
\hline 
\end{tabular}
\end{center}
\end{table}

With Table II  in hand, it is now possible to increase the parameter $a$ to expand the range of validity. 
Several examples of this are given in Table III.

\enlargethispage{40pt}
\begin{table}[!h]
\caption{More derived and very widely applicable effective bounds on $\pi(x)$.\\
\centerline{(Here instead $a$ is adjusted so as to maximize the region of validity, $x_0=2$.)}\\
\centerline{$|\pi(x)-\Li(x)| <   a  \;x \;(\ln x)^{ b} \; \exp\left(-{ c}\; \sqrt{\ln x}\right).$}
}\smallskip
\begin{center}
\begin{tabular}{||c|c|c|c||}
\hline
\hline
$a$ & $b$ & $c$ & $x_0$  \\
\hline
\hline
 3/4& --3/4 &0.32115 & 2 \\
 \hline
 7/8& --7/8 &0.32115 & 2 \\
\hline
 0.93& --1 &0.32115 & 2 \\
\hline
\hline
0.8935& --3/4 &0.3936 & 2  \\ 
\hline
0.94 & --7/8 &0.3936 & 2  \\ 
\hline
1.05& --1 &0.3936 & 2 \\
\hline
\hline
\end{tabular}
\end{center}
\end{table}

We shall subsequently apply these bounds (both Tables II and III) in various ways --- specifically to bounding the location of the $n^{th}$ prime. 

\section{Exponentially bounding the prime number theorem}

Note that for any $a>0$, $b\geq -1$, $c>0$, we can write:
\begin{equation}
 a \;x \;(\ln x)^b \; \exp\left(-c\; \sqrt{\ln x}\right) 
 = a\; {x\over\ln x}\; (\ln x)^{b+1} \; \exp\left(-c\; \sqrt{\ln x}\right).
\end{equation}
It is useful to define
\begin{equation}
f(x) = (\ln x)^{b+1} \; \exp\left(-c\; \sqrt{\ln x}\right).
\end{equation}

Note  that  ${x\over\ln x} < \pi(x)$ for  $x>17$, so that
\begin{equation}
 a \;x \;(\ln x)^b \; \exp\left(-c\; \sqrt{\ln x}\right) < a\; \pi(x) \; f(x); \qquad (x>17).
 \end{equation}
 
 \enlargethispage{30pt}
 This is already somewhat interesting since it implies that the various bounds on $|\pi(x)-\Li(x)|$ discussed above all lead to 
 \begin{equation}
|\pi(x)-\Li(x)| <  a  \; \pi(x) \;f(x); 
\qquad (x \geq \max\{x_0,17\}).
\label{E:pi-Li}
\end{equation}
 
 When feasible, ($x_0$ sufficiently small), direct computation might potentially increase the domain of validity
\begin{equation}
|\pi(x)-\Li(x)| <  a  \; \pi(x) \;f(x); 
\qquad (x \geq x_*; \; x_*\leq\max\{x_0,17\}).
\end{equation}
This can be rearranged to yield
\begin{equation}
{\Li(x) \over 1+a\; f(x)} < \pi(x) < {\Li(x) \over 1-a \;f(x)}; \qquad (x \geq x_*).
\end{equation}
Note that $f(x)\to 0$ with exponential rapidity in $\sqrt{\ln x}$. 

To be more explicit, for $x \geq x_*$:
 \begin{equation}
{\Li(x) \over 1+a\; (\ln x)^{b+1} \; \exp\left(-c\; \sqrt{\ln x}\right)} < \pi(x) < {\Li(x) \over 1-a \;(\ln x)^{b+1} \; \exp\left(-c\; \sqrt{\ln x}\right)}.
\end{equation}
These inequalities  provide a remarkably clean and effective version of the prime number theorem. 

\clearpage
\section{Localizing the $n^{th}$ prime:}

Note that for $b\geq -1 $ the function $f(x) = (\ln x)^{b+1} \; \exp(-c\; \sqrt{\ln x})$
achieves a global maximum at $x_{peak} = \exp\left( \left[2(b+1)/ c\right]^2\right)\geq 1$, above which $f'(x)<0$, and below which $f'(x)>0$. 

Concentrating on the region above the peak we shall first re-write $f(x)$ in terms of $p_{\pi(x)}$, yielding a bound in terms of  $f(p_{\pi(x)})$.
We will then convert this into an explicit function of $\pi(x)$, yielding a bound in terms of $f( \pi(x) \ln \pi(x))$. 
Doing so will require some careful attention to the domain of validity of the various inequalities.
(One wishes the various manipulations to bring one closer to the peak without overshooting the peak.)

Observe that  in all generality $p_{\pi(x)} \leq x < p_{\pi(x)+1}$,  and more specifically one has 
$p_{\pi(x_{peak})} \leq x_{peak} < p_{\pi(x_{peak})+1}$
So if we now choose $x > p_{\pi(x_{peak})+1} > x_{peak}$ then we can safely write 
\begin{equation}
f(x) < f(p_{\pi(x)}); \qquad ( x > p_{\pi(x_{peak})+1} \}).
\end{equation}
In 1952 Nagura~\cite{Nagura} (among other results) showed
\begin{equation}
\theta\left({4\over 3} \; x\right) -\theta(x) >0; \qquad (x\geq 109 = p_{29}).
\end{equation}
This bound is by no means optimal, but it will be good enough for current purposes.
This implies
\begin{equation}
p_{n+1} < {4\over 3} \; p_n; \qquad (n\geq 29).
\end{equation}
Explicitly checking smaller integers yields
\begin{equation}
p_{n+1} < {4\over 3} \; p_n; \qquad (n\geq 5=\pi(11)).
\end{equation}
But then let us choose $x>{4\over 3} \; x_{peak}$ and note
\begin{equation}
x>{4\over 3} x_{peak} \geq {4\over 3} p_{\pi(x_{peak})} > p_{\pi(x_{peak})+1}
\end{equation}
Then certainly
\begin{equation}
f(x) < f(p_{\pi(x)}); \qquad \left( x > \max\left\{ 11, \; {4\over3} \;x_{peak} \right\}\right).
\end{equation}

\enlargethispage{40pt}
Now invoke Rosser's theorem $p_n > n \ln n$, ($n \geq 1$)~\cite{Rosser}, then 
\begin{equation}
\label{E:Rosser}
f(x) <  f(\pi(x)\ln \pi(x)); 
\end{equation}
where in addition to the prior constraints  $x > \max\{ 11, \; {4\over3} x_{peak}\}$ we now need to check that $\pi(x)\ln(\pi(x))> x_{peak}$.

To get a more computationally useful grasp on the domain of validity use the fact that 
$p_n < n \ln(n \ln n)$ for $n\geq 6$~\cite{Rosser2} and note
\begin{equation}
p_n < n \ln(n \ln n) = n(\ln n + \ln\ln n) = n \ln n \left\{1 + {\ln\ln n\over\ln n} \right\};
\qquad (n \geq 6).
\end{equation}
The quantity in braces is bounded and takes on the maximum value $1+e^{-1}$ at $x=e^e$.
Thence
\begin{equation}
p_n < (1+e^{-1}) \;  n \ln n = 1.367879441 \;  n \ln n;
\qquad (n \geq 4).
\end{equation}
So we see
\begin{equation}
\pi(x)\ln(\pi(x)) > {p_{\pi(x)}\over 1+ e^{-1}};  \qquad (x\geq p_4 = 7)
\end{equation}
This implies that  the condition $\pi(x)\ln(\pi(x))> x_{peak}$ relevant to inequality (\ref{E:Rosser}) is certainly satisfied for $ p_{\pi(x)} > (1+e^{-1}) x_{peak}$, which in turn is certainly satisfied for $ x > (1+e^{-1}) x_{peak}> {4\over 3} x_{peak}$, so that we can safely write
\begin{equation}
f(x) <  f(\pi(x)\ln \pi(x)); \qquad (x > \max\{11, \;  (1+e^{-1}) \; x_{peak}\}).
\end{equation}
Finally, combine this result with the conditions for the validity of the bound (\ref{E:pi-Li}), \\to yield 
\begin{equation}
|\pi(x)-\Li(x)| <  a\; \pi(x) \; f(\pi(x)\ln \pi(x)); \qquad (x\geq \max\{ x_0, 17,  (1+e^{-1}) x_{peak}\}).
\label{E:final?}
\end{equation}
We remind the reader that for current purposes we are interested in
\begin{equation}
x_{peak} = \exp\left( \left[2(b+1)\over c\right]^2\right) \geq 1.
\end{equation}
Because we have made various approximations getting to this stage, the actual range of the validity of the final inequality (\ref{E:final?}) may be somewhat larger than naively advertised $\max\{ x_0, 17, (1+e^{-1}) x_{peak}\}$. 
That is, while the bounds (\ref{E:final?}) are certainly guaranteed to hold for values of $x$ sufficiently far above the peak, they might still hold both at the peak and for some region below the peak.

\enlargethispage{40pt}
Whenever feasible (meaning whenever $x_{peak}$ is not too large) this should be checked by explicit computation to find a suitable $x_*$ such that:
\begin{equation}
|\pi(x)-\Li(x)| <  a\; \pi(x) \; f(\pi(x)\ln \pi(x)); \qquad 
(x\geq x_*; \; x_* \leq \max\{ x_0, 17,  (1+e^{-1}) x_{peak}\}).
\label{E:final!}
\end{equation}

Sometimes $x_*$ is significantly lower than the naive bound $\max\{ x_0, 17, (1+e^{-1})  x_{peak}\}$, and
sometimes $ (1+e^{-1})  x_{peak}\leq 17$ and so that the peak can be safely ignored. But sometimes we simply need to keep track of this complication.

If we evaluate this final inequality (\ref{E:final!}) at $x=p_n$ then
\begin{equation}
|\Li(p_n)-n| <  a\; n \; f(n\ln n); \qquad (n\geq n_*; \; n_* \leq \max\{ \pi(x_0), 7, \pi( (1+e^{-1})\;x_{peak})\}).
\end{equation}
Finally, as promised in the abstract,
\begin{equation}
 \Li^{-1}(  n [1-a \; f(n\ln n)])< p_n < \Li^{-1}(  n [1+a \; f(n\ln n)]); \qquad (n\geq n_*). 
 \end{equation}
 To be  more explicit
 \[
p_n < \Li^{-1} \left( n \left[1+  a   \;(\ln[n\ln n])^{b+1} \; \exp\left(-{c}\; \sqrt{\ln[n\ln n]}\right)\right] \right);
\qquad (n\geq n_*);
\]
\[
p_n > \Li^{-1} \left( n \left[1-  a   \;(\ln[n\ln n])^{b+1} \; \exp\left(-{c}\; \sqrt{\ln[n\ln n]}\right)\right] \right);
\qquad (n\geq n_*).
\]
Here the range of validity is explicitly bounded by some $n_*$ with 
\[
n_* \leq \max\left\{\pi(x_0),7, \pi\left( (1+e^{-1}) \exp\left( \left[2(b+1)\over c\right]^2\right)\right) \right\}.
\]
These effective bounds can be viewed as fully explicit versions of the usual Cippola asymptotic expansion --- see for instance the discussion in reference~\cite{asymptotically}.

See Table IV below (based on  Tables I, II, and III) for some explicit examples.
The outputs from Table II, $\{\tilde a,\tilde b,\tilde c, \tilde x_0\}$, have been relabelled as $\{a,b,c,x_0\}$, and used as inputs for new computations of $x_*$ and $n_*$ for the new inequalities discussed in  this section.

\vspace{-10pt}
\enlargethispage{75pt}
\begin{table}[!h]
\caption{Some widely applicable effective bounds on $|\Li(x)-\pi(x)|$ and $p_n$:\\
\centerline{$|\Li(x)-\pi(x)| < a\; \pi(x) \; \;(\ln[\pi(x) \ln \pi(x) ])^{b+1} \; \exp\left(-{c}\; \sqrt{\ln[\pi(x)\ln \pi(x)]}\right);$} \\
\centerline{$p_n \lessgtr \Li^{-1} \left( n \left[1\pm  a   \;(\ln[n\ln n])^{b+1} \; \exp\left(-{c}\; \sqrt{\ln[n\ln n]}\right)\right] \right)$. }}\smallskip
\begin{center}
\begin{tabular}{||c|c|c||c|c|c||c|c|c||}
\hline
\hline
$a$ & $b$ & $c$ & $x_0$ & $x_{peak}$ & $x_*$ & $n_0$  & $n_{peak}$ & $n_*$\\
\hline
\hline
0.4394& --3/4 &0.32115 & 59 & 11.3 & 41 &17  & 5 & 13 \\ 
\hline
 0.4680 & --7/8 & 1/4  & 41& $e$ & 37 & 13 & 1& 12 \\ 
\hline
 0.4795 & --1 & 1/6  & 41 & 1 &37  &  13 & 0 &12 \\
\hline
\hline
0.2795& --3/4 &0.3936 & 229 & 5.022 &  149 &50 & 3 & 35 \\ 
\hline
0.2804 &--5/6 & 1/3  & 227 & $e$ & 149 & 49 & 1 &35\\
\hline
 0.3164  & --1 & 1/4  & 223 & 1 & 97 & 48 & 0 & 25 \\
\hline
\hline
9.2211 & 1/2 & 0.8476 & 2 & 275 789& 2 & 1 & 24 104 & 1\\ %
\hline
 9.7590 & 0 & 1/2& 2 & $e^{16}$ &2 &1 & 595 341 & 1 \\ 
 \hline
11.9026 & --1/2 & 1/5  & 2 & $e^{25}$ & * & 1 &\vphantom{\Big | } ${1\over25}e^{25}$ &*\\ %
\hline 
29.6698 & --1 & 1/10  & 2 &1 &2 & 1  & 0 & 1 \\ %
\hline
\hline
9.59 & 0.515 & 0.8274 & 2 &667 161 & 2 & 1  & 54 105 & 1\\
\hline
11.148 & 0 & 1/2 & 2 & $e^{16}$ & 2 & 1 &595 431& 1\\
\hline
13.659 & --1/2 & 1/5 & 2 &  $e^{25}$ &  *& 1 &\vphantom{\Big | } ${1\over25}e^{25}$ &*\\ %
\hline
34.955 & --1 & 1/10 & 2 &  1& 2 & 1 & 0 & 1\\
\hline
\hline 
\end{tabular}\\
Explicit evaluation of $x_*$ and $n_*$ in the situations labelled $*$ proved \\
computationally infeasible due to the magnitude of $x_{peak}$ and $n_{peak}$.\\
All indications are that $x_*=2$ and $n_*=1$ in both cases.
\end{center}
\end{table}
 
 \clearpage
\section{Special case}\label{S:special-case}

As previously indicated, the special case $b=-1$ is particularly appealing. 

Starting from any bound of the form
\begin{equation}
|\pi(x)-\Li(x)| <  a  \;x \;(\ln x)^{-1} \; \exp\left(-{c}\; \sqrt{\ln x}\right); 
\qquad (x \geq x_0);
\end{equation}
we have the particularly simple result  that:
\begin{equation}
|\pi(x)-\Li(x)| <  a  \; \pi(x)\; \exp\left(-{c}\; \sqrt{\ln x}\right); 
\qquad (x \geq x_*; \; x_* \leq \max\{x_0,17\}).
\end{equation}
Thence in the range $x\geq x_*$ we have
 \begin{equation}
{\Li(x) \over 1+a\; \exp\left(-c\; \sqrt{\ln x}\right)} < \pi(x) < {\Li(x) \over 1-a \; \exp\left(-c\; \sqrt{\ln x}\right)}.
\end{equation}
Furthermore, (noting that in this case we always have $x_{peak}=1$),
\begin{equation}
p_n < \Li^{-1} \left( n \left[1+  a  \exp\left(-{c}\; \sqrt{\ln[n\ln n]}\right)\right] \right);
\qquad (n\geq n_*);
\end{equation}
\begin{equation}
p_n > \Li^{-1} \left( n \left[1-  a   \exp\left(-{c}\; \sqrt{\ln[n\ln n]}\right)\right] \right);
\qquad (n\geq n_*),
\end{equation}
subject to  the simple constraint
\begin{equation}
n_*  \leq \max\{\pi(x_0),7\}.
\end{equation}
For several explicit effective examples of this specific behaviour see Table V (which is a subset of Table IV and a minor extension of parts of Table~III).
These bounds are not necessarily optimal but, (given the initial input information in Table I), they are relatively simple to derive.

\enlargethispage{40pt}
\begin{table}[!h]
\caption{Widely applicable effective bounds on $|\Li(x)-\pi(x)|$ and $p_n$ for $b=-1$:\\
\centerline{$|\Li(x)-\pi(x)| < a\; \pi(x) \; \; \exp\left(-{c}\; \sqrt{\ln[\pi(x)\ln \pi(x)]}\right);$} \\
\centerline{$p_n \lessgtr \Li^{-1} \left( n \left[1\pm  a   \; \exp\left(-{c}\; \sqrt{\ln[n\ln n]}\right)\right] \right)$. }}\smallskip
\begin{center}
\begin{tabular}{||c|c|c||c|c|c||c|c|c||}
\hline
\hline
$a$ & $b$ & $c$ & $x_0$ & $x_{peak}$ & $x_*$ & $n_0$  & $n_{peak}$ & $n_*$\\
\hline
\hline
 0.4795 & --1 & 1/6  & 41 & 1 &37  &  13 & 0 &12 \\
\hline
\hline
 0.3164  & --1 & 1/4  & 223 & 1 & 97 & 48 & 0 & 25 \\
\hline
\hline
29.6698 & --1 & 1/10  & 2 &1 &2 & 1  & 0 & 1 \\ %
\hline
\hline
34.955 & --1 & 1/10 & 2 &  1& 2 & 1 & 0 & 1\\
\hline
\hline 
\hline
\hline
 1& --1 &0.32115 & 2 &1& 2 & 1 & 0 & 1 \\
\hline
\hline
1.1 & --1 &0.3936 & 2 & 1& 2 & 1 & 0 & 1\\
\hline
\hline
\end{tabular}
\end{center}
\end{table}
 
\section{Further developments}\label{S:further}

Another  useful trick is to evaluate the inequality
\begin{equation}
|\pi(x)-\Li(x)| <  a\; \pi(x) \; f(\pi(x)\ln \pi(x)); \qquad (x\geq x_*)
\end{equation}
both at $x=p_n$ and at $x=p_{n+1}^-$, just below $p_{n+1}$. 
Then since $\pi(p_{n+1}^-)=n=\pi(p_n)$ we have both
\begin{equation}
|n-\Li(p_n)| <  a\; n \; f(n\ln n); \qquad (n\geq n_*);
\end{equation}
and
\begin{equation}
|n-\Li(p_{n+1})| <  a\; n \; f(n\ln n); \qquad (n\geq n_*).
\end{equation}
But then by the triangle inequality
\begin{equation}
\Li(p_{n+1})-\Li(p_n) <  2 a\; n \; f(n\ln n); \qquad (n\geq n_*);
\end{equation}
Implying
\begin{equation}
p_{n+1} < \Li^{-1} \left(\Li(p_n) +  2 a\; n \; f(n\ln n)\right); \qquad (n\geq n_*). 
\end{equation}
This observation can be used to develop yet more asymptotic bounds along the lines of the Cippola expansion, as discussed in~\cite{asymptotically}, and can also be related to the effective bounds developed using the first Chebyshev $\theta(x)$ function as described in reference~\cite{prime-gaps}.

\section{Conclusions}\label{S:discussion}

With some hindsight, an ineffective version of the bounds considered herein could have been developed some 125 years ago. The current analysis is interesting in two respects: First, using developments made over the past decade, the bounds have been made explicitly and fully effective. Secondly, because the bounds are asymptotically exponential, they will eventually always overtake and outperform bounds based on the Cippola expansion.
Finally, I do not claim any of these bounds are in any way optimal. As always there is a three-way trade-off between stringency of the bound, simplicity of the bound, and the region of validity. Herein I have largely focussed on simplicity of the bound and a wide (sometimes maximal) region of validity.

\clearpage

\end{document}